\documentclass[a4paper,12pt]{amsart}
\usepackage{amssymb}
\usepackage{ifthen}
\usepackage{graphicx}
\usepackage{float}
\usepackage{caption}
\usepackage{subcaption}
\usepackage{cite}
\usepackage{amsfonts}
\usepackage{amscd}
\usepackage{amsxtra}
\usepackage{mathrsfs}
\usepackage[usenames]{color}

\newtheorem{thm}{Theorem}
\newtheorem{cor}{Corollary}
\newtheorem{lem}{Lemma}

\newtheorem{rem}{Remark}
\newtheorem{example}{Example}
\newtheorem{defn}{Definition}
\newtheorem{prob}{Problem}

\newtheorem{conj}{Conjecture}
\theoremstyle{definition}

\newcounter {own}
\def\theown {\thesection  .\arabic{own}}

\newenvironment{pf}[1][]{%
 \vskip 3mm
 \noindent
 \ifthenelse{\equal{#1}{}}%
  {{\slshape Proof. }}%
  {{\slshape #1.} }%
 }%
{\qed\bigskip}

\newcounter{alphabet}
\newcounter{tmp}
\newenvironment{Thm}[1][]{\refstepcounter{alphabet}%
\bigskip%
\noindent%
{\bf Theorem \Alph{alphabet}}%
\ifthenelse{\equal{#1}{}}{}{ (#1)}%
{\bf .} \itshape}{\vskip 8pt}

\makeatletter
\newcommand{\Ref}[1]{\@ifundefined{r@#1}{}{\setcounter{tmp}{\ref{#1}}\Alph{tmp}}}
\makeatother

\newenvironment{Lem}[1][]{\refstepcounter{alphabet}%
\bigskip%
\noindent%
{\bf Lemma \Alph{alphabet}}%
{\bf .} \itshape}{\vskip 8pt}

\newcommand{\IR}{{\mathbb R}}
\newcommand{\IN}{{\mathbb N}}

\newcommand{\ID}{{\mathbb D}}

\def\be{\begin{equation}}
\def\ee{\end{equation}}

\newcommand{\bee}{\begin{enumerate}}
\newcommand{\eee}{\end{enumerate}}

\newcommand{\blem}{\begin{lem}}
\newcommand{\elem}{\end{lem}}
\newcommand{\bthm}{\begin{thm}}
\newcommand{\ethm}{\end{thm}}
\newcommand{\bcor}{\begin{cor}}
\newcommand{\ecor}{\end{cor}}
\newcommand{\beg}{\begin{example}}
\newcommand{\eeg}{\end{example}}
\newcommand{\begs}{\begin{examples}}
\newcommand{\eegs}{\end{examples}}
\newcommand{\bdefe}{\begin{defn}}
\newcommand{\edefe}{\end{defn}}
\newcommand{\bprob}{\begin{prob}}
\newcommand{\eprob}{\end{prob}}
\newcommand{\bei}{\begin{itemize}}
\newcommand{\eei}{\end{itemize}}

\newcommand{\bcon}{\begin{conj}}
\newcommand{\econ}{\end{conj}}
\newcommand{\bcons}{\begin{conjs}}
\newcommand{\econs}{\end{conjs}}
\newcommand{\bprop}{\begin{propo}}
\newcommand{\eprop}{\end{propo}}
\newcommand{\br}{\begin{rem}}
\newcommand{\er}{\end{rem}}
\newcommand{\brs}{\begin{rems}}
\newcommand{\ers}{\end{rems}}
\newcommand{\bo}{\begin{obser}}
\newcommand{\eo}{\end{obser}}
\newcommand{\bos}{\begin{obsers}}
\newcommand{\eos}{\end{obsers}}
\newcommand{\bpf}{\begin{pf}}
\newcommand{\epf}{\end{pf}}
\newcommand{\ba}{\begin{array}}
\newcommand{\ea}{\end{array}}
\newcommand{\beq}{\begin{eqnarray}}
\newcommand{\beqq}{\begin{eqnarray*}}
\newcommand{\eeq}{\end{eqnarray}}
\newcommand{\eeqq}{\end{eqnarray*}}

\newcommand{\ds}{\displaystyle}

\newcounter{minutes}
\divide\time by 60
\newcounter{hours}
\multiply\time by 60 \addtocounter{minutes}{-\time}

\begin{document}
\bibliographystyle{amsplain}
\title[]{Sections of univalent harmonic mappings}

\thanks{
File:~\jobname .tex,
          printed: \number\year-\number\month-\number\day,
          \thehours.\ifnum\theminutes<10{0}\fi\theminutes}

\author{Saminathan Ponnusamy
}
\address{S. Ponnusamy, and A. Sairam Kaliraj,
Indian Statistical Institute (ISI), Chennai Centre, SETS, MGR Knowledge City, CIT
Campus, Taramani, Chennai 600 113, India. }
\email{samy@isichennai.res.in, sairamkaliraj@gmail.com}
\author{Anbareeswaran Sairam Kaliraj}

\author{Victor V. Starkov}
\address{V.V. Starkov, Petrozavodsk State University, 33, Lenin Str., 185910, Petrozavodsk, Republic of Karelia, Russia.}
\email{vstarv@list.ru}


\date{\today  
}

\begin{abstract}
In this article,
we determine the radius of univalence of sections of normalized univalent harmonic mappings for which the range is convex (resp. starlike, close-to-convex, convex in one direction). Our result on the radius of univalence of section $s_{n,n}(f)$ is sharp
especially when the corresponding mappings have convex range.
In this case, each section $s_{n,n}(f)$ is univalent in the
disk of radius $1/4$ for all $n\geq2$, which may be compared with classical result of Szeg\"{o} on conformal mappings.

\end{abstract}
\thanks{ }

\maketitle
\pagestyle{myheadings}
\markboth{S.Ponnusamy, A. Sairam Kaliraj, and V.V. Starkov}{Sections of univalent harmonic mappings}
\section{Introduction and main results}\label{PS7Sec1}
Since confirmation of the Bieberbach conjecture by Louis de Branges \cite{de_Branges} on the class
$\mathcal{S}$, of all normalized univalent analytic functions $\phi$ defined in the unit disk $\mathbb{D} = \{z: \; |z|<1\}$, one of the
open problems about the class $\mathcal S$ is  that of determining
the precise value of $r_n$ such that all sections $s_n(\phi)$ of $\phi$ are univalent in $|z|<r_n$. Here
we say that $\phi$ is normalized if $\phi(0)=0=\phi'(0)-1$. Also, let
\begin{equation}\label{p7_an_sec}
s_n(\phi)(z)= \sum_{k=1}^{n} a_kz^k
\end{equation}
whenever
\begin{equation}\label{ser_rep_ana}
\phi(z)= \sum_{k=1}^{\infty} a_kz^k.
\end{equation}
In \cite{Szego}, Szeg\"{o} proved that the section/partial sum $s_n(\phi)$ of $\phi\in\mathcal{S}$ is univalent in $|z| < 1/4$ for all $n \geq 2$.
The constant $1/4$ is sharp as the second section of the Koebe function $k(z)=z/(1-z)^2$ suggests.
In \cite{Rober41}, Robertson proved that the section $s_n(k)$ is starlike in the disk $|z|<1-3n^{-1}\log n$ for $n \geq 5$, and that the number $3$ cannot be replaced by a smaller constant.
Later in the year 1991, Bshouty and Hengartner \cite{Bshouty-Hengartner} showed that the Koebe function is not extremal for the problem of determining the radius of univalency of the partial
sums of functions in $\mathcal{S}$. At this time the best known result is due to Jenkins \cite{Jenkins} who proved that $s_n(\phi)$ is univalent in $|z|<r_n$ for
$\phi\in\mathcal{S}$, where the radius of univalence $r_n$ is at least $ 1-(4\log n - \log(4\log n))/n$ for $n \geq 8$. For related investigations on this topic,
see the recent articles \cite{ObSamy13,Hiroshi-Samy-2014} and the references therein.
More interestingly, as investigated recently in \cite{LiSamyNA1, LiSamyNA2}, our main aim in this article is to consider the analogous
problem for univalent harmonic mappings in the unit disk, since harmonic mappings have interesting links with geometric
function theory, minimal surfaces and locally quasiconformal mappings.

Every harmonic mapping $f$ in a simply connected domain can be written as $f=h+\overline{g}$, where  $h$ and $g$ are analytic. In particular,
we consider the class ${\mathcal H}$ of all complex-valued harmonic functions $f=h+\overline{g}$ in
${\mathbb D}$ normalized by $h(0)=g(0)=0=h'(0)-1 $. We call $h$ and $g$, the analytic and the co-analytic
parts of $f$, respectively, and obviously they have the following power series representation
\be\label{PSSerRep}
h(z)=z+\sum _{k=2}^{\infty}a_k z^k ~\mbox{ and }~ g(z)=\sum _{k=1}^{\infty}b_k z^k, ~z \in \mathbb{D}.
\ee
Throughout the discussion we shall use this representation.
Since the Jacobian $J_f$ of $f=h+\overline{g}$ is $J_f(z) = |h'(z)|^2-|g'(z)|^2,$
we say that $f$ is sense-preserving in $\ID$ if $J_f(z)>0$ in $\ID$. Let $\mathcal{S}_H$ denote the class of all sense-preserving
harmonic univalent mappings $f \in {\mathcal H}$ and set $\mathcal{S}^0_H=\{f \in \mathcal{S}_H:\,  f_{\overline{z}}(0)=0\}$.
For many basic results on univalent harmonic mappings, we refer to the monograph of  Duren \cite{Duren:Harmonic} and also \cite{PonRasi2013}. Harmonic mappings techniques have been used to study and solve fluid flow problems (see \cite{Aleman_Constantin}). In particular, the study of univalent harmonic functions having special geometric property such as convexity, starlikeness, and close-to-convexity  arises naturally while dealing with planar fluid dynamics problems. For example, in \cite[Theorem 4.5]{Aleman_Constantin},  Aleman et al. considered a fluid flow problem on a convex domain $\Omega_0$ satisfying an interesting geometric property. In view of results from \cite{PonSai-1(11),Po-Hiroshi-Hirosh-2013}, one obtains that harmonic mappings,
$z\mapsto H_0(z) +\overline{G_0(z)}$ such that ${\rm Re}\,H_0'(z)>|G_0'(z) |$ on the convex domain $\Omega_0$,
considered by the authors in \cite[Theorem 4.5]{Aleman_Constantin} are indeed close-to-convex in $\Omega_0$.

Another reasons in studying the sections of harmonic mappings is that
approximation of real valued harmonic functions by harmonic polynomials attracted the attention of mathematicians (see \cite{Walsh})
as it has many advantages. For example, a harmonic function has its maximum and minimum values on the boundary of the regions
of consideration.  Because planar harmonic mappings $f=h+\overline{g}$ defined on $\ID$ have series representation as in \eqref{PSSerRep},
sections of $f$ can be thought of as an approximation of $f$ by complex-valued harmonic polynomials and thus,
approximation of univalent harmonic mappings by univalent harmonic polynomials might lead to new applications
in fluid flow problems, in particular. Until recently, much is not known about the univalence of sections of univalent harmonic mappings.
In 2013, Li and Ponnusamy \cite{LiSamyNA1, LiSamyNA2} initiated the study on this topic by considering certain classes of univalent harmonic
mappings. However, the harmonic analog of these results are not known in the literature even for well-known geometric
subclasses of $\mathcal{S}^0_H$, namely, the classes ${\mathcal S}_H^0$ are ${\mathcal K}_H^0$, ${\mathcal S}_H^{*0}$, and ${\mathcal C}_H^0$
mapping $\ID$ onto, respectively, convex, starlike, and close-to-convex domains, just as
${\mathcal K}$, ${\mathcal S}^*$, and ${\mathcal C}$ are the subclasses of ${\mathcal S}$
mapping $\ID$ onto these respective domains. At this place, it is worth to recall that  general theorems on convolutions \cite{Rusch_Sheil} (see also \cite[p.~256, 273]{Duren}) allow one to conclude that $s_n(\phi)$ is convex, starlike, or close-to-convex in the disk $|z|<1-3n^{-1}\log n$, for $n \geq 5$, whenever $\phi\in {\mathcal S}$ is convex (resp. starlike or close-to-convex) in $\ID$.

Another interesting geometric subclass of ${\mathcal C}_H^0$
which attracted function theorists is the class of univalent harmonic functions $f$ for which $f(\ID)$
is convex in a direction $\alpha$. Recall that a domain $D\subset\mathbb{C}$ is called convex
in the direction $\alpha$ $(0\leq \alpha< \pi)$ if the intersection of $D$ with each line
parallel to the line through $0$ and $e^{i\alpha}$ is connected (or empty).
See, for example, \cite{Clunie-Small-84,HengSch70,RZ}. Now, we recall the class $\mathcal{S}^0_H(\mathcal{S})$ introduced in \cite{PonSai5}, where
$$\mathcal{S}^0_H(\mathcal{S}) = \left\{h+\overline{g} \in \mathcal{S}^0_H :\,
h+e^{i \theta}g \in \mathcal{S}~ \mbox{for some}~~\theta \in \mathbb{R} \right\}$$
and as in \cite{PonSai5}, let
$\mathcal{S}_H(\mathcal{S})=\left\{f=f_0 +b \overline{f_0}:\,f_0 \in \mathcal{S}^0_H(\mathcal{S})~~\mbox{and}~~b \in \mathbb{D} \right\}.
$
One of the conjectures stated in \cite{PonSai5} reads as follows.

\begin{conj}\label{PS5conj2}
$\mathcal{S}^0_H = \mathcal{S}^0_H(\mathcal{S})$. That is, for every function $f=h+\overline{g}
\in \mathcal{S}^0_H$, there exists at least one $\theta \in \mathbb{R}$ such that $h + e^{i\theta} g \in \mathcal{S}$.
\end{conj}

In \cite{PonSai5}, it was also remarked that the truth of this conjecture verifies the
coefficient conjecture of Clunie and Sheil-Small for $f=h+\overline{g}\in {\mathcal S}_H^0$, namely,
$$|a_n|\le\frac {(2n+1)(n+1)} {6}, \quad |a_{-n}|\le\frac {(2n-1)(n-1)} {6}, \quad \mbox{and} \quad
 \big ||a_n|-|a_{-n}|\big |\le n
$$
for each $n\geq 2$, where $b_n=a_{-n}$.  Conjecture \ref{PS5conj2} remains open. The bound $|a_{-2}|\le 1/2$  is well-known and sharp
which follows from the classical Schwarz lemma. However,
the conjectured bounds of  Clunie and Sheil-Small have been verified for a number of subclasses of  ${\mathcal S}_H^0$,
namely, for ${\mathcal S}_H^{*0}$,  ${\mathcal C}_H^0$,
$\mathcal{S}^0_H(\mathcal{S})$ and the class of harmonic mappings convex in one direction.
More recently, Starkov \cite{Starkov} established a criteria for
functions belonging to the class  $\mathcal{S}^0_H(\mathcal{S})$
and as a consequence, several examples including harmonic univalent polynomials are also obtained for a given $f\in \mathcal{S}_H$.

For $f=h+\overline{g} \in \mathcal{S}^0_H$ with power series representation as in \eqref{PSSerRep}, the sections/partial sums
$s_{n,m}(f)$ of $f$ are defined as
$$s_{n,m}(f)(z)=s_n(h)(z)+\overline{s_m(g)(z)}
$$
where $n\geq 1$ and $m\geq 2$. However, the special case $m=n\geq 2$ seems interesting in its own merit. We now state our main results.

\bthm\label{PS7Thm2}
Let $f=h+\overline{g} \in \mathcal{S}^0_H$ with series representation as in \eqref{PSSerRep}. Suppose that $f$ belongs to any one of the
following geometric subclasses of ${\mathcal S}_H^0$ : ${\mathcal S}_H^{*0}$,  ${\mathcal C}_H^0$,
$\mathcal{S}^0_H(\mathcal{S})$ or the class of harmonic mappings convex in one direction.
Then the section $s_{n,m}(f)$ is univalent in the disk $|z|<r_{n,m}$. Here $r_{n,m}$ is the unique positive
root of the equation $\psi(n,m,r)=0$, where
\be\label{PS7_eq1}
\psi(n,m,r) = \frac {1}{12 r}
\left(\frac {1-r}{1+r}\right)^3 \left[1-\left(\frac {1-r}{1+r}\right)^{6}\right] - R_n - T_m,
\ee
with
\be\label{PS7_th2_eq1.1}
R_n = \sum_{k=n+1}^{\infty}A_k r^{k-1}, ~~ T_m = -\sum_{k=m+1}^{\infty}A_{-k} r^{k-1}, ~\mbox{ where }~
A_k=\frac{k(k+1)(2k+1)}{6},
\ee

In particular, every section $s_{n,n}(f)(z)$ is univalent in the disk $|z|<r_{n,n}$, where
$$r_{n,n} > r^L_{n,n}:=1- \frac{(7\log n  - 4\log(\log n))}{n} ~\mbox{ for }~ n \geq 15.
$$
Moreover, $r_{n,m} \geq r^L_{l,l}$,  where  $l=\min\{n,m\}\geq 15$.
\ethm

For functions in the convex family $\mathcal{K}^0_H$ of harmonic mappings, we have the following interesting result which may compared with
the original conjecture for functions in $\mathcal S$.

\bthm\label{PS7Thm3}
Let $f=h+\overline{g} \in \mathcal{K}^0_H$ with series representation as in \eqref{PSSerRep}. Then the section $s_{n,m}(f)$ is univalent in the disk $|z|<r_{n,m}$, where $r_{n,m}$ is the unique positive root of the equation $\mu(n,m,r)=0$. Here
\be\label{PS7_eq5}
\mu(n,m,r)=\frac {1-r}{(1+r)^3} - \sum_{k=n+1}^{\infty}\left[ \frac{k(k+1)}{2} r^{k-1} \right] - \sum_{k=m+1}^{\infty}\left[ \frac{k(k-1)}{2} r^{k-1} \right].
\ee
In particular, for $n\geq5$, and $\theta \in \IR$, the harmonic function
$$s_{n,n}(f;\theta)(z)=s_n(h)(z)+e^{i\theta}\overline{s_n(g)(z)}
$$
is univalent and close-to-convex in the disk $|z| < 1-3n^{-1}\log n $.
Moreover, we have $r_{n,m} \geq 1-(4\log l - 2\log(\log l))/l$, where $l=\min\{n,m\} \geq 7 $.
\ethm

It is worth to remark that if $f \in \mathcal{K}^0_H$, then we actually prove that for $n\geq5$, $s_{n,n}(f)$ is stable harmonic
close-to-convex (see \cite{Rodri-Maria}) in $|z|<1-3n^{-1}\log n$.

The paper is organized as follows. In Section \ref{PS7Sec2}, we recall certain known results which are crucial in the proof of our main theorems.
In Section \ref{PS7Sec4}, we present the proofs of Theorems \ref{PS7Thm2} and \ref{PS7Thm3}, and as a consequence, we state a  couple of corollaries.

\section{Useful Lemmas}\label{PS7Sec2}

Now, we recall some results that are needed for the proofs of our main results.
The following result due to Bazilevich \cite{Bazilevich} gives the necessary and sufficient condition for a normalized analytic function to be univalent in $\ID$.

\begin{Thm}\label{uni_nec_suf_Analytic}
An analytic function $\phi$ defined in $\ID$ and determined by \eqref{ser_rep_ana} is univalent in $\ID$ if and only if for each
$z \in \ID$ and each $t\in[0, \pi/2]$,
\be\label{PS7inteq1}
\frac {\phi(re^{i\eta})-\phi(re^{i\psi})} {re^{i\eta}-re^{i\psi}} :=\sum_{k=1}^{\infty}a_k \frac{\sin kt}{\sin t} z^{k-1} \ne 0,
\ee
where $t=(\eta-\psi)/2$, $z=re^{i(\eta+\psi)/2}$ and $\left .\frac{\sin kt}{\sin t}\right |_{t=0}=k$.
\end{Thm}

Recently,   Starkov  \cite{Starkov} generalized this result to the class of normalized sense-preserving harmonic mappings in the following form.

\begin{Thm}\label{uni_nec_suf}
A sense-preserving harmonic function $f=h+\overline{g}$ defined in $\ID$ determined by \eqref{PSSerRep} is univalent in $\ID$ if and only if for each
$z \in \ID\setminus\{0\}$ and each $t\in(0, \pi/2]$,
\be\label{PS7inteq2}
 \frac {f(re^{i\eta})-f(re^{i\psi})} {re^{i\eta}-re^{i\psi}}  :=\sum_{k=1}^{\infty}\left[ (a_k z^k - \overline{b_k z^k})\frac{\sin kt}{\sin t} \right] \ne 0,
\ee
where $t=(\eta-\psi)/2$ and $z=re^{i(\eta+\psi)/2}$.
\end{Thm}


%

The following two point distortion theorem of Graf et al. \cite{Graf-samy}
plays a crucial role in the proof of our main results.

\begin{Lem}\label{two_po_dis}
If $f=h+\overline{g} \in \mathcal{S}^0_H$, $r \in (0, 1)$, $t, \psi \in \IR$, then
$$\left|\frac {f(re^{it})-f(re^{i\psi})} {re^{it}-re^{i\psi}}\right|\ge \frac {1}{4\alpha r}
\left(\frac {1-r}{1+r}\right)^\alpha \left[1-\left(\frac {1-r}{1+r}\right)^{2\alpha}\right],
$$
where $\alpha = \sup_{f \in \mathcal{S}_H}|h''(0)|/2$.
\end{Lem}

Finally, we  recall the following well-known identities which are also easy to derive.
\blem\label{idenities1}
The following identities are true for $0 < r < 1$:
\begin{itemize}
\item[(i)] $\ds \sum _{k=n+1}^{\infty}k r^{k-1} = \frac{r^n}{(1-r)^2}[1+n(1-r)]$.
\item[(ii)] $\ds \sum _{k=n+1}^{\infty}k^2r^{k-1} = \frac{r^n}{(1-r)^3}[2+(2n-1)(1-r)+n^2(1-r)^2]$.
\item[(iii)] $\ds \sum _{k=n+1}^{\infty}k^3r^{k-1} = \frac{r^n}{(1-r)^4}[6+(6n-6)(1-r)+(3n^2 -3n +1)(1-r)^2+n^3(1-r)^3]$.
\end{itemize}
\elem
\bpf
The identity (i) is obvious. To obtain (ii) we may multiply (i) by $r$ and then differentiate it with respect to $r$.
The proof of case (iii) is similar. So, we omit its proof.
\epf

\section{Partial sums of Univalent Harmonic Mappings}\label{PS7Sec4}
\subsection{Proof of Theorem \ref{PS7Thm2}}
Suppose that $f=h+\overline{g}$ belongs to either ${\mathcal S}_H^{*0}$ or ${\mathcal C}_H^0$ or
$\mathcal{S}^0_H(\mathcal{S})$ or to the class of harmonic mappings convex in one direction, where  $h, g$ are given by the
power series \eqref{PSSerRep} with $b_1=0$. Set $F_r(z)=f(rz)/r$ for $0 < r < 1$.
Then  $F_r(z) \in \mathcal{S}^0_H$ and
$$ F_r(z) = z+\sum _{k=2}^{\infty}a_k r^{k-1} z^k + \sum _{k=2}^{\infty}\overline{b_k r^{k-1} z^k}.
$$
Evidently, finding the largest radius of univalence of $s_{n,m}(f)(z)$ is equivalent to finding the largest value $r$
such that $s_{n,m}(F_r)(z)$ is univalent in $\ID$. From Theorem \Ref{uni_nec_suf}, it is clear that $s_{n,m}(F_r)(z)$
is univalent in $\ID$ if and only if $s_{n,m}(F_r)(z)$ is sense-preserving in $\ID$ and the associated section $P_{n,m,r}(z)$
has the property that
$$P_{n,m,r}(z) := \sum_{k=1}^{M}\left[ (a'_k z^k - \overline{b'_k z^k})\frac{\sin kt}{\sin t} \right] \ne 0, ~\mbox{ for all }~ z\in\ID\setminus\{0\}, ~\mbox{ and }~
t\in(0,\pi/2],
$$
where $M = \max\{n, m\}$, $l=\min\{n, m\}$, $a'_k = a_k r^{k-1}$, $b'_k = b_k r^{k-1}$ for all $k \leq l$,
$$
a'_k= \left \{ \begin{array}{lr}
\ds a_k r^{k-1} &  \mbox{for all $k > l$ if $M = n$},\\
0 & \mbox{for all $k > l$ if $M > n$},
\end{array}
\right.
$$
and
$$
b'_k=\left \{ \begin{array}{lr}
\ds b_k r^{k-1} &  \mbox{for all $k > l$ if $M = m$},\\
0 & \mbox{for all $k > l$ if $M > m$}.
\end{array}
\right.
$$
Setting $t=(\eta-\psi)/2$, $z=\rho e^{i(\eta+\psi)/2}\in\ID$ in \eqref{PS7inteq2} and from the univalency of $F_r$ for $0<r < 1$, we get that
\be\label{PS7_eq3}
\left|\sum_{k=1}^{\infty}\left[ (a_k z^k - \overline{b_k z^k})r^{k-1}\frac{\sin kt}{\sin t} \right]\right| \geq  \frac {1}{12 r}
\left(\frac {1-r}{1+r}\right)^3 \left[1-\left(\frac {1-r}{1+r}\right)^{6}\right].
\ee
In order to find a lower bound for $|P_{n,m,r}(z)|$, we need to find an upper bound for
$$\left|R_{n,m,r}(z)\right| = \left| \sum_{k=n+1}^{\infty}\left[ a_k r^{k-1} z^k \frac{\sin kt}{\sin t} \right] - \sum_{k=m+1}^{\infty}\left[\overline{(b_k r^{k-1} z^k)} \frac{\sin kt}{\sin t} \right]\right|.
$$
By the assumption on $f$,  it follows that (see for instance \cite{Duren:Harmonic,PonSai5})
$$
|a_k| \leq \frac{(k+1)(2k+1)}{6} ~\mbox{ and }~ |b_k| \leq \frac{(k-1)(2k-1)}{6}, ~\mbox{ for all }~ k \geq 2,
$$
and hence
\beq\label{PS7_eq3.1}\nonumber
|R_{n,m,r}(z)| &\leq& \sum_{k=n+1}^{\infty}\left[ \frac{k(k+1)(2k+1)}{6} r^{k-1} \right] + \sum_{k=m+1}^{\infty}\left[ \frac{k(k-1)(2k-1)}{6} r^{k-1} \right] \\
 &=& R_n + T_m \qquad \mbox{(see \eqref{PS7_th2_eq1.1})}.
\eeq
From \eqref{PS7_eq3} and \eqref{PS7_eq3.1}, we get that
$$|P_{n,m,r}(z)| \geq \frac {1}{12r} \left(\frac {1-r}{1+r}\right)^3 \left[1-\left(\frac {1-r}{1+r}\right)^{6}\right] - R_n - T_m =\psi(n,m,r).
$$
The inequality $|P_{n,m,r}(z)| > 0$ holds for all $z\in\ID\setminus\{0\}$,  whenever $\psi(n,m,r)>0$, where $\psi(n,m,r)$ is
defined by \eqref{PS7_eq1}. This gives that  $\psi(n,m,r)>0$ for all $r \in (0, r_{n,m})$, where
$r_{n,m}$ is the positive root of the equation $\psi(n,m,r)=0$ which lies in the interval $(0, 1)$.
In order to complete the proof, we have to show that $s_{n,m}(f)$ is locally univalent in $|z|< r_{n,m}$. However,
$s_{n,m}(f)=s_n(h)+\overline{s_m(g)}$ is locally univalent in $|z|< r_{n,m}$ if and only if the analytic functions $s_n(h)+e^{i\theta}s_m(g)$ is locally univalent in $|z|< r_{n,m}$ for every $\theta \in \IR$. That is, we have to show that $s_{n-1}(h')(z)+e^{i\theta}s_{m-1}(g')(z) \ne 0$ for all $|z|< r_{n,m}$ and $\theta \in \IR$.
It is easy to see that
$$s_{n-1}(h')(z)+e^{i\theta}s_{m-1}(g')(z)= (h'(z)+e^{i\theta}g'(z)) - \left(\sum_{k=n+1}^{\infty} k a_k z^{k-1} + e^{i\theta} \sum_{k=m+1}^{\infty} k b_k z^{k-1}\right).
$$
From the hypothesis, it is clear that $f=h+\overline{g}$ belongs to either $\mathcal{C}^0_H$ or $\mathcal{S}^0_H(S)$ (see for instance, \cite{PonSai5}).
As the affine spanning of $\mathcal{C}^0_H$ as well as $\mathcal{S}^0_H(S)$ are linear invariant families \cite{PonSai5},
$$\alpha:=\sup_{f\in \mathcal{C}_H \bigcup \mathcal{S}_H(S)}\frac{|h''(0)|}{2} = 3.
$$
Therefore, from a well known result on linear invariant family of harmonic mappings,  it follows that
(see \cite[p. 99]{Duren:Harmonic})
$$|h'(z)| - |g'(z)| \geq \frac{(1-r)^2}{(1+r)^4} ~\mbox{ for $|z|=r$ }.
$$
Moreover, for $0 < r < 1$ (see \eqref{PS7_eq3}), we see that
$$\min \left\{\frac{(1-r)^2}{(1+r)^4}, ~\frac {1}{12 r}
\left(\frac {1-r}{1+r}\right)^3 \left[1-\left(\frac {1-r}{1+r}\right)^{6}\right] \right\} = \frac {1}{12 r}
\left(\frac {1-r}{1+r}\right)^3 \left[1-\left(\frac {1-r}{1+r}\right)^{6}\right]
$$
and thus,
\beqq
\left|s_{n-1}(h')(z)+e^{i\theta}s_{m-1}(g')(z)\right| &\geq&  \frac{(1-r)^2}{(1+r)^4} - R_n - T_m \\
&\geq& ~\frac {1}{12 r}
\left(\frac {1-r}{1+r}\right)^3 \left[1-\left(\frac {1-r}{1+r}\right)^{6}\right]- R_n - T_m \\
&\geq& 0,
\eeqq
where the last inequality here gives the condition $|z|<r_{n,m}$. This observation proves that $s_{n,m}(f)$ is univalent in the disk $|z|< r_{n,m}$
and the proof of the first part of the theorem is complete.

From the above discussion, it is apparent that $r_{n,m} \geq r_{l,l}$, where $l=\min\{n, m\}\geq 2$.
Next, we need to consider the special case $m=n$ and determine the lower bound for $r_{n,n}$ with certain restriction on $n$.
In this case, the sufficient condition \eqref{PS7_eq1} for the univalence of $s_{n,n}(f)$ reduces to
\be\label{PS7_eq4}
\psi(n,n,r)=\frac{(1-r)^3 (3 + 10 r^2 + 3 r^4)}{3(1+r)^9} - (R_n + T_n) \geq 0,
\ee
where $R_n + T_n$ determined from \eqref{PS7_th2_eq1.1} may be rearranged in a convenient form as
\beqq
R_n + T_n &=& \frac{1}{3}\sum_{k=n+1}^{\infty} k r^{k-1} + \frac{2}{3} \sum_{k=n+1}^{\infty} k^3 r^{k-1}\\
&=& \frac{r^n [12+12(n-1)(1 - r)+3(2n^2-2n+1)(1 - r)^2+(2n^3+n)(1 - r)^3]}{3 (1 - r)^4}.
\eeqq
As $\lim_{n\rightarrow\infty}(R_n + T_n) =0$, it is clear that the radius of univalence $r_{n,n}\rightarrow 1$.
Setting $r = 1-\alpha_n/n$, where $\alpha_n=o(n)$, we see that \eqref{PS7_eq1} holds, whenever $0<t(\alpha_n, n)<1$,
where
$$t(x, n)=\frac{e^{-x}n^7}{x^7}\left(2 - \frac{x}{n}\right)^9\frac{[12n^4+12(n-1)xn^3+3(2n^2-2n+1)x^2n^2+(2n^3+n)x^3n]}{16 n^4 -
32 n^3 x + 28 n^2 x^2 - 12 n x^3 + 3 x^4}.
$$
From the definition of $t(\alpha_n, n)$, it is clear that $e^{\alpha_n}t(\alpha_n, n)=O(n^7)$ and hence $\alpha_n$ can be chosen to be
$7\log n  - a\log(\log n)$ for some positive real number $a$. However, computations shows that
$$\lim_{n\rightarrow\infty} t(\alpha_n, n) = \frac{64}{2401} > 0 ~\mbox{for}~ a = 4,
~\mbox{ and }~
\lim_{n\rightarrow\infty} t(\alpha_n, n) = \infty ~\mbox{for}~ a > 4.
$$
Therefore, we may set $\gamma_n=7\log n  - 4\log(\log n).$

It is easy to see that $1- \gamma_n/n > 0$ for all $n \geq 15$. For $n \geq 15$, we shall prove that $r_{n,n}>1- \gamma_n/n$.
For $n \geq 15$, it is sufficient to prove that $t(x, n)$ is a decreasing function in $x$, whenever
$$\gamma_n \leq  x \leq n,~~0 < t(\gamma_n, n) < 1 ~\mbox{ and }~ t(n, n)>0.
$$
In order to do this, we first differentiate $t(x,n)$ with respect to $x$ and obtain that $t'(x, n)=q_1(x,n) q_2(x,n)$, where
\beqq
q_1(x,n) &=& \frac{-(2 n - x)^8 e^{-x} }{x^8[16 n^4 -  32 n^3 x + 28 n^2 x^2 - 12 n x^3 + 3 x^4]^2}
\eeqq
and
\beqq
q_2(x,n) &=& 2688 n^7 + 2688(n-3)n^6x + 3(448 n^7-2368n^6+3648n^5-x^7)x^2\\
&~~& + 64n^4(7n^3-48 n^2+137n-132)x^3+16 n^2 (59 n^3 - 128 n^2 + 178 n -75 ) x^5\\
&~~& + 2 n^2(32n^5-80n^4(6+x)+1672n^3-4n^2(774+13 x^3)+2040n-3x^5)x^4\\
&~~& +2n(88n^4-240n^3+434n^2-390n+81) x^6\\
&~~& +2n(78n^2-98n+57)x^7+ 6(6n^3-2n^2+6n-1)x^8.
\eeqq
From the definition of $q_1(x,n)$, it is clear that $q_1(x,n) < 0 $ for all $\gamma_n \leq  x \leq n$, where $n \geq 15$.
To conclude $t'(x, n) < 0$, we need to prove that $q_2(x,n)>0$ for all $n \geq 15$ and $\gamma_n \leq  x \leq n$.
From the grouping of terms in $q_2(x,n)$, one can see that $q_2(x,n)>0$ for $n \geq 15$, and $x \in (0, n/3]$.
Next, we show that $q_2(x,n)>0$ for all $x \in [n/3, n]$. To do this, for any fixed $n\geq15$, we set $x=n/k$, where $k\in[1, 3]$.
Then, $q_2(x,n)$ reduces to
$$Q(k, n)=\frac{n^7(1-2 k)^2}{k^6}Q_1(k)+\frac{n^8}{k^8}Q_2(k)+\frac{n^9}{k^9}Q_3(k)+\frac{n^{10}}{k^8}Q_4(k)+\frac{n^{11}}{k^9}Q_5(k),
$$
where
\beqq
Q_1(k)&=& 6(112k^4-224k^3+204k^2-92k+27), \\
Q_2(k)&=& 2(1344k^7-3552k^6+4384k^5-3096k^4+1424k^3-390k^2+57k-3),\\
Q_3(k)&=& 1344 k^7-3072k^6+3344k^5-2048k^4+868k^3-196k^2+36k-3, \\
Q_4(k)&=& 4(112 k^5-240k^4+236k^3-120k^2+39k+3)~~\mbox{ and }\\
Q_5(k)&=& 2(32k^5-80k^4+88k^3-52k^2+18k-3).
\eeqq
A computation shows that $Q_1(k)$ has no real root. Moreover, the only real roots of $Q_2(k)$, $Q_3(k)$, $Q_4(k)$ and $Q_5(k)$ are $0.104153$,
$0.143187$, $-0.0630667$ and $0.5$, respectively. Therefore, for $1\leq j \leq 5$, $Q_j(k)$ will have same sign for all $k\in[1, 3]$.
As $Q_j(1)>0$ for $1\leq j\leq 5$, we conclude that $Q_j(k)>0$ for all $k\in[1, 3]$. This shows that $Q(k, n)>0$ for all $k \in [1, 3]$ and
hence, $q_2(x,n)>0$ for all $x \in [n/3, n]$. Since
$$t(n, n)=\frac{e^{-n}}{3}(2 n^3+ 6 n^2+ 7 n +3) > 0 \mbox{ for all }~ n\in \IN,
$$
from the fact that $t'(x, n)<0$, we infer that $t(x, n)$ is a positive and decreasing function of $x$ in the interval $(0, n]$, for each $n\geq 15$.
To complete the proof, we have to show that $t(\gamma_n, n)<1$ for all $n\geq 15$. By making use of upper bounds of $\gamma_n/n$
for various values of $n$, it is easy to see that $t(\gamma_n, n)<1$ for $n\geq 73$. A direct computation using mathematica shows
that $t(\gamma_n, n)<1$ for $15 \leq n < 73$ also. The proof is complete.
\hfill$\Box$

\bcor
Let $f \in \mathcal{S}^0_H$ satisfies the hypothesis of Theorem \ref{PS7Thm2}. Then $s_{n,n}(f)(z)$ is univalent in the disk
\begin{itemize}
\item[(i)] $|z|<1/4$, whenever $n \geq 7$,
\item[(ii)] $|z|<1/2$, whenever $n \geq 22$,
\item[(ii)] $|z|<3/4$, whenever $n \geq 78$.
\end{itemize}
The bound for the radius of univalence $r_{n,n}$ of $s_{n,n}(f)$ for certain values of $n$ are listed in Table \ref{tab2}.
\ecor
\begin{table}
\begin{center}
\begin{tabular}{|l|l|}
  \hline
  Value of $n$ & Value of $r_{n,n}$ \\
  \hline
  2 & 0.108193 \\
  \hline
  3 & 0.147197 \\
  \hline
  4 & 0.182263 \\
  \hline
  5 & 0.214025 \\
  \hline
  10 & 0.337088 \\
  \hline
  50 & 0.675001 \\
  \hline
  100 & 0.788521 \\
  \hline
  287 & 0.900122 \\
  \hline
\end{tabular}
\end{center}
\caption{Values of $r_{n,n}$ for certain values of $n$\label{tab2}}
\end{table}

%

The following shearing theorem due to Clunie and Sheil-Small is needed for the proof of Theorem \ref{PS7Thm3}.

\begin{Thm}\label{lem2.1}{\rm \cite[Theorem 5.3]{Clunie-Small-84}}
A locally univalent harmonic function $f=h+\overline{g} $  in $\mathbb{D}$ is a
univalent mapping of $\mathbb{D}$ onto a domain convex in the
direction $\theta$ if and only if $\phi_\theta=h-e^{i2\theta}g$ is a conformal
univalent mapping of $\mathbb{D}$ onto a domain convex  in the
direction $\theta$.
\end{Thm}

The proof of Theorem \ref{PS7Thm3} is similar to the proof in Theorem \ref{PS7Thm2} with the help of the corresponding
coefficients inequalities and the sharp lower bound for the two point distortion theorem (Lemma \ref{two_po_dis_har_con}) for $f \in \mathcal{K}^0_H$.


\begin{lem}\label{two_po_dis_har_con}
If $f = h+\overline{g} \in \mathcal{K}^0_H$, $r \in (0, 1)$, $t, \psi \in \IR$, then
$$\left|\frac {f(re^{it})-f(re^{i\psi})} {re^{it}-re^{i\psi}}\right|\ge \frac {1-r}{(1+r)^3}.
$$
\end{lem}
\bpf
For every pair of points $re^{it}$ and $re^{i\psi}$ in $\ID$, we can find a $\theta \in \IR$ such that
$$(h(re^{it})-h(re^{i\psi}))+\overline{(g(re^{it})-h(re^{i\psi}))} = (h(re^{it})-h(re^{i\psi}))+e^{i\theta}(g(re^{it})-h(re^{i\psi})).
$$
Since $f \in \mathcal{K}^0_H$, $f$ is convex in every direction and hence, by Lemma \Ref{lem2.1}, the function $h-e^{i2\theta}g$
is univalent in $\ID$ for every $\theta \in \IR$. The desired conclusion follows from the two point distortion theorem
for univalent analytic functions (see \cite[Corollary 7, p.~127]{Duren}).
\epf

\subsection{Proof of Theorem \ref{PS7Thm3}}
Let $f=h+\overline{g} \in \mathcal{K}^0_H$. Then the Taylor coefficients of $h$ and $g$ satisfy the inequality (see \cite{Clunie-Small-84,Duren:Harmonic})
\be\label{PS7_thm3_eq_1}
|a_n|\leq\frac{n+1}{2} ~\mbox{ and }~ |b_n|\leq\frac{n-1}{2} ~\mbox{ for all }~ n \geq 2.
\ee
Following the proof of Theorem \ref{PS7Thm2} with the same notation and \eqref{PS7_thm3_eq_1}, we get that
$$|R_{n,m,r}(z)| \leq  \sum_{k=n+1}^{\infty}\left[ \frac{k(k+1)}{2} r^{k-1} \right] + \sum_{k=m+1}^{\infty}\left[ \frac{k(k-1)}{2} r^{k-1} \right] = R_n + T_m.
$$
From Lemma \ref{two_po_dis_har_con}, we obtain that
$$|P_{n,m,r}(z)| \geq \frac {1-r}{(1+r)^3} - R_n - T_m.
$$
The inequality $|P_{n,m,r}(z)|>0$ holds for all $z\in \ID\setminus\{0\}$, whenever $\mu(n,m,r)>0$, where $\mu(n,m,r)$ is defined by \eqref{PS7_eq5}.
However, $\mu(n,m,r)>0$ for all $r\in(0, r_{n,m})$, where $r_{n,m}$ is the unique positive root of the equation \eqref{PS7_eq5} which is less
than $1$. Since the affine span of $\mathcal{K}^0_H$ is a linear invariant family and $\alpha = \sup_{f \in \mathcal{K}_H}|h''(0)|/2 = 2$, we have (see \cite[p. 99]{Duren:Harmonic})
$$|h'(z)| - |g'(z)| \geq \frac{(1-r)}{(1+r)^3} ~\mbox{ for }~ |z|=r.
$$
Following the proof technique of Theorem \ref{PS7Thm2}, we conclude that $s_{n,m}(f)$ is locally univalent in $|z| < r_{n,m}$. Now, let us first consider the special case $m=n$. In this case, \eqref{PS7_eq5} reduces to
\beqq
\mu(n,n,r)=\frac {1-r}{(1+r)^3} - \sum_{k=n+1}^{\infty} k^2 r^{k-1}.
\eeqq
From Lemma \ref{idenities1}, the expression for $\mu(n,n,r)$ simplifies to
\be\label{PS7th1eq2}
\mu(n,n,r)= \frac{1-r}{(1+r)^3} - \frac{ r^n [2+(2n-1)(1-r)+n^2(1-r)^2]}{(1-r)^3}.
\ee
For $0 < r < 1$, $\mu(n,n,r) > 0$ if and only if $ \psi(n, r) > 0 $, where
$$\psi(n, r) = (1-r)^4 - [2+(2n-1)(1-r)+n^2(1-r)^2] (1+r)^3 r^n. $$ From the continuity of $\psi(n, r)$ and from the fact that $\psi(n, 0) > 0$ and $\psi(n, 1) < 0$ , it is evident that there exists a real number $r_n>0$ such that $\psi(n, r) > 0 $ for all $r \in (0, r_n)$ and $\psi(n, r_n)=0$. To complete the proof, we need to find the lower bound for $r_n$ for large values of $n$. By letting $r=1-\alpha_n/n$ in $\psi(n, r)$ and making use of the fact that $e^{-\alpha_n/n}\geq1-\alpha_n/n$, for $\alpha_n \geq 0$, we see that $\psi(n, r) > 0$ holds whenever $0 < T(\alpha_n, n) < 1$,
where
\be\label{PS7th1eq1_1}
T(x, n)= \frac{e^{-x} n^4}{{x}^4}\left(2-\frac{x}{n}\right)^3\left(2 + (2n-1)\frac{x}{n} + {x}^2\right).
\ee
Clearly $\alpha_n=o(n)$ and $e^{\alpha_n}T(\alpha_n, n)=O(n^4)$ and hence the dominating term in $\alpha_n$ could be
at most $4\log n$. This allows us to choose $\alpha_n = 4\log n - a\log\log n$ for appropriate choices of $a>0$. As $\lim_{n \rightarrow \infty} T(\alpha_n, n) = \infty$ for $a > 2$, $a$ cannot be chosen to be  greater than $2$ in our estimate. Moreover, for $a = 2$, a computation shows that $\lim_{n \rightarrow \infty} T(\alpha_n, n) = 1/2 < 1$. Therefore, we may set
$$\beta_n = 4\log n - 2\log(\log n) ~\mbox{ for }~ n\geq7.
$$
Throughout the further discussion in this proof, we assume that $n\geq7$ since $1-\beta_n/n > 0$ holds only when $n\geq7$.

For $\beta_n \leq  x \leq n$, we shall prove that $0<T(x, n)<1$. We observe that it suffices to show that
$T'(x, n) < 0$  for $\beta_n \leq  x \leq n$, $0 < T(\beta_n, n) < 1$  and $0 < T(n, n) < 1$.
Indeed, differentiating \eqref{PS7th1eq1_1} with respect to $x$ one can obtain by a standard calculation that
$$T'(x, n) = \frac{- (2 n - x)^2 [2 n^2 (8 + 8 x + 4 x^2 + x^3)-n x(8 + 4 x + x^2 + x^3) + x^3]}{e^{x}x^5} <0,
$$
whenever $\beta_n < x \leq n$. As $(2-\beta_n/n)>0$ for all $n \geq 7$, we deduce that
$$T(\beta_n, n)=\left( \frac{\log n}{{\beta_n}^2}\right)^2\left(2-\frac{\beta_n}{n}\right)^3\left[2 + (2n-1)\frac{\beta_n}{n}
+ {\beta_n}^2\right] > 0 ~\mbox{ for all }~ n \geq 7.
$$
Next, we show that $T(\beta_n, n) < 1$ for all $n \geq 7$. To do this, we may  rewrite the expression for $T(\beta_n, n)=T_1(\beta_n, n)+ T_2(\beta_n, n)+ T_3(\beta_n, n)$,
where
\beqq
T_1(\beta_n, n) &=&  \frac{16(\log n)^2}{{\beta_n}^4}\left(1-\frac{\beta_n}{2n}\right)^3,\\
T_2(\beta_n, n) &=&  \frac{16(\log n)^2}{{\beta_n}^3}\left(1-\frac{\beta_n}{2n}\right)^3\left(1-\frac{1}{2n}\right),~\mbox{ and }~\\
T_3(\beta_n, n) &=&  \frac{8(\log n)^2}{{\beta_n}^2}\left(1-\frac{\beta_n}{2n}\right)^3.
\eeqq
First, we show that  $T(\beta_n, n) <  1$ for all $n\geq 16$.
We introduce
\beqq
A(x)&=&\log x -\log(\log x)+ \frac{(\log(\log x))^2}{4\log x}\\
B(x)&=&\log x - \frac{\log(\log x)}{2}, ~\mbox{ and }~C(x)=2 - \frac{\log(\log x)}{\log x}.
\eeqq
We see that $A(x)$ is an increasing function of $x$ in the interval $[7, \infty)$ and thus, we
easily have
$$A(n)=\log n -\log(\log n)+ \frac{(\log(\log n))^2}{4\log n} \geq  A(9) >\sqrt{2} ~\mbox{ for }~ n\geq9.
$$
Similarly, it is easy to see that $B(x)$ is increasing in the interval $[7, \infty)$ and $C(x)$
is increasing in the interval $[16, \infty)$. Therefore, a simple computation shows that
$$B(n)=\log n - \frac{\log(\log n)}{2} \geq B(16)\approx 2.2627 ~\mbox{ for}~ n\geq 16
$$
and
$$C(n)=2-\frac{\log(\log n)}{\log n} \geq C(16) \approx 1.63219 ~\mbox{ for}~ n\geq 16.
$$
By making use of the above inequalities, let us find an upper bound for each of the quantities
$T_1(\beta_n, n)$, $T_2(\beta_n, n)$ and $T_3(\beta_n, n)$. We begin with
\beqq
T_1(\beta_n, n) &=&  \frac{16(\log n)^2}{{\beta_n}^4}\left(1-\frac{\beta_n}{2n}\right)^3 <  \frac{1}{{\beta_n}^2[1-\frac{\log(\log n)}{2\log n}]^2}\\
&=&  \frac{1}{16[\log n -\log(\log n)+ \frac{(\log(\log n))^2}{4\log n}]^2} < \frac{1}{32} ~~\mbox{for all}~ n\geq 9.
\eeqq
Next,  we obtain
\beqq
T_2(\beta_n, n) &=&  \frac{16(\log n)^2}{{\beta_n}^3}\left(1-\frac{\beta_n}{2n}\right)^3\left(1-\frac{1}{2n}\right)\\
&=&  \frac{1}{[\log n - \frac{\log(\log n)}{2}][2 - \frac{\log(\log n)}{\log n}]^2} < \frac{1}{6} ~~\mbox{for all}~ n\geq 16.
\eeqq
Finally, for all $n\geq 16$,
\beqq
T_3(\beta_n, n) =  \frac{8(\log n)^2}{{\beta_n}^2}\left(1-\frac{\beta_n}{2n}\right)^3
<  \frac{2}{[2-\frac{\log(\log n)}{\log n}]^2} < \frac{19}{24} .
\eeqq
Using these bounds, it clear that for all $n\geq 16$ one has $T(\beta_n, n) <1$.
On the other hand, a direct computation yields that $T(\beta_n, n) < 1 $ for $7\leq n \leq 15$.
Hence $s_{n,n}(f)$ is univalent in $|z| < r_{n,n}$, where $r_{n,n} >  1-(4\log n - 2\log(\log n))/n$ for $n \geq 7$.
From the above proof technique it is apparent that $r_{n,m} \geq r_{l,l}$, where $l=\min\{n, m\}$. That is, for $l=\min\{n,m\} \geq 7 $, we have $r_{n,m} \geq 1-(4\log l - 2\log(\log l))/l$.

However, geometric properties of convex harmonic mappings gives a way to improve the lower bounds of $r_{n,n}$. From Theorem \Ref{lem2.1},
it is clear that $\phi_\theta(z)=h(z)-e^{i2\theta}g(z)$ is univalent and convex in the direction of $\theta$.
As univalent functions convex in one direction are close-to-convex in $\ID$, it is clear that $\phi_\theta$ is close-to-convex in $\ID$
for every $\theta \in \IR$. From a classical result on convolution of analytic functions \cite[Theorem 8.7, p. 248]{Duren},
it is clear that the radius of close-to-convexity of $s_n(\phi_\theta)$ is greater than or equal
to the radius of convexity of $s_n(l(z))$, where $l(z)=z/(1-z)$. From \cite[Corollary 3, p.~256]{Duren},
it is clear that the radius of univalence of $s_n(\phi_\theta)$ cannot exceed the radius of convexity of $s_n(l(z))$.
The radius of convexity of $s_n(l(z))$ is known to be $1-3n^{-1}\log n $ for $n\geq5$.

But then, from a result of Clunie and Sheil-Small \cite[Lemma 5.15]{Clunie-Small-84}, we get that
$$s_{n,n}(f;\theta)(z)=s_{n,n}(h)(z)+e^{i\theta}\overline{s_{n,n}(g)(z)}
$$
is univalent and close-to-convex in the disk $|z|<1-3n^{-1}\log n $ for all $n\geq5$.
\hfill $\Box$

\vspace{8pt}


\bcor
Let $f \in \mathcal{K}^0_H$ satisfy the hypothesis of Theorem \ref{PS7Thm3}. Then $s_{n,n}(f;\theta)(z)$ is univalent in the disk
\begin{itemize}
\item[(i)] $|z|<1/4$, whenever $n \geq 2$,
\item[(ii)] $|z|<1/2$, whenever $n \geq 17$,
\item[(ii)] $|z|<3/4$, whenever $n \geq 46$.
\end{itemize}
\ecor
\bpf
The proof of Case (i) follows from Lemma \Ref{lem2.1} and the result of Szeg\"{o} \cite{Szego}. The rest of the cases follows from Theorem \ref{PS7Thm3}.
\epf
\subsection*{Acknowledgements}
 The research was supported by the project RUS/RFBR/P-163 under
Department of Science \& Technology (India) and the Russian Foundation for Basic Research (project 14-01-92692). The first
author is currently on leave from Indian Institute of Technology Madras.

\end{document}